\documentclass[11pt,reqno]{amsart}

\usepackage{latexsym}
\usepackage{amsfonts}
\usepackage{amsmath}
\usepackage{url}
\usepackage{graphicx}
\usepackage{color}
\usepackage{bbm}
\definecolor{myblue}{rgb}{0.09,0.32,0.44} 
\usepackage[pdfborder={0 0 0},pdfborderstyle={},colorlinks=true,linkcolor=myblue,citecolor=myblue,urlcolor=blue]{hyperref}

\usepackage{epsfig,psfrag,verbatim,amsfonts}

\usepackage{amssymb,amsmath}

\def\bl{\begin{lemma}}
\def\el{\end{lemma}}
\def\bth{\begin{theorem}}
\def\eth{\end{theorem}}
\def\bc{\begin{corollary}}
\def\ec{\end{corollary}}
\def\bcj{\begin{conjecture}}
\def\ecj{\end{conjecture}}
\def\bpr{\begin{proposition}}
\def\epr{\end{proposition}}
\def\bde{\begin{definition}}
\def\ede{\end{definition}}

\newcommand{\be}{\begin{eqnarray}}
\newcommand{\ee}{\end{eqnarray}}

\newcommand{\Z}{{\mathbb Z}}

\newcommand{\SL}{{\mathrm{SL}}}
\newcommand{\Cay}{{\mathrm{Cay}}}

\renewcommand{\and}{\hbox{ {\rm and} }}

\newtheorem{theorem}{Theorem}[section]
\newtheorem{definition}{Definition}[section]
\newtheorem{lemma}[theorem]{Lemma}

\newtheorem{corollary}[theorem]{Corollary}
\newtheorem{proposition}[theorem]{Proposition}
\newtheorem{conjecture}[theorem]{Conjecture}

\newtheorem*{theorem*}{Theorem}

\theoremstyle{definition}
\numberwithin{equation}{section}

\begin{document}
\title{A note on the structure of expanders}
\author{Itai Benjamini and Mikolaj Fraczyk}

\date{20.10.20}

\maketitle

\begin{abstract}
A conjecture regarding the structure of expander graphs is discussed. 
\end{abstract}

\section{}

 A finite graph $G$, is called $h$-{\it expander} for
$$ 
h = \inf_{S \subset V_G :  \mbox{ } 0 < |S| < |G|/2} {|\partial S|
\over |S|},
$$

where $V_{G}$ are the vertices of $G$ and $\partial S$ is the outer vertex boundary of $S$.

 $G$ a graph with degrees bounded from above by $d$ and expansion $h$.

\begin{conjecture} There is a function $f(h, d) >0$ so that $G$ contains a spanning subgraph, which is an $f(h,d)$- expander with girth proportional to the diameter. Where the proportionality ratio between girth and diameter,
are bounded below by strictly positive functions only of the degrees and the expansion of $G$.

\end{conjecture}

Start by showing existence of spanning expander with  girth going to infinity with the size of $G$? 
\medskip

Note that we don't require from the spanning high girth sub expander to be regular. 
\medskip

We don't know the conjecture even further assuming that the expanders are Cayley graphs. 
\medskip

Are there approachable high dimensional formulations of the conjecture for the several variants of high dimensional expanders \cite{L}?
\medskip 

For the standard definitions and background see \cite{HLW}.

\section{Motivation and examples}

\begin{enumerate}

	\item 
	In \cite{BS} it was shown that infinite graphs with expansion bigger or equal $1$, admit a spanning tree with with expansion bigger or equal $1$. This is an infinite variant on the conjecture for sufficiently large expansion. For finite expanders with expansion at least $1$, it implies that there are balls of order diameter with expanding spanning trees. Maybe using the Lovasz local lemma one can glue these to get a spanning subexpander?
	\medskip
	
	\item
	One could try a random construction, trimming short cycles.

	Under some further assumptions already Bernoulli percolation gives a giant component which is rather tree like locally. By adapting theorem 4 of \cite{pyond} one can conclude that if $p$ is such that 
	$$
	\rho(G)d_{G}p < 1.
	$$
	
Where $\rho(G)$ is the spectral radius, $d_{G}$ the degree of $G$. And $p$-Bernoulli bond percolation admits a giant component. Then there will be a unique giant component (\cite{abs})
	which locally , up to order diameter, has the structure of an infinite percolation cluster in the non-uniqueness regime, admitting  dense triforcation vertices. For any expander $G$, if one takes a sufficiently large power of $G$ the percolation requirements  will hold \cite{PS}.
	\medskip
	
	Show that if $G^2$ satisfies the conjecture, then $G$
	as well.
	
	For trimming of percolation clusters to get a subgraph with positive expansion in the context of infinite vertex transitive graphs see \cite{BLS}.   
	\medskip
	
		\item
	
One specific  challenge is products of 
graphs  $G \times G$, were $G$
	is an expander, e.g. a random regular graph or a Ramanujan graph.
	\medskip

	As a very partial positive indication note that
	
 $\SL(m,Z/p^nZ)\times \SL(m,Z/p^nZ)$, $m>1$ where $p$ is prime and $n$ goes to infinity admits Cayley graphs with girth going to infinity.
Using Corollary 1.4 from \cite{BG} we find a free dense, finitely generated subgroup $\Lambda \in \SL(m,\mathbb Z_p)\times \SL(m,\mathbb Z_p)$. Let $S\subset \Lambda$ be a finite set that freely generates $\Lambda$.

$\Lambda$ is dense, so for each $n$ the image $S$ generates $\SL(m,\mathbb Z/p^n\mathbb Z)\times \SL(m,\mathbb Z/p^n\mathbb Z)$. Consider a sequence of Cayley graphs 
$$
Cay(\SL(m,\mathbb Z/p^n\mathbb Z)\times \SL(m,\mathbb Z/p^n\mathbb Z, S)).
$$

The girth of this sequence tends to infinity. As otherwise there would be a finite length relation on elements of $S$ that holds modulo $p^n$ for every $n$. But then the same relation would have to hold in the inverse limit which is $SL(m,\mathbb Z_p)\times SL(m,\mathbb Z_p)$. This contradicts the freeness.
	\medskip

	\item Let $S\subset \SL(2,\mathbb Z)$ be a finite symmetric set generating a Zariski dense subgroup. By the work of Bourgain and Gamburd \cite{BoGa} it is known that the sequence of graphs $G_p:=\Cay(\SL(2,\mathbb Z/p\mathbb Z), S)$ is an expander sequence as $p$ varies among sufficiently big primes. Let $\Gamma$ be the subgroup generated by $S$. If $S$ is free then the argument form the previous example shows that $G_p$ have large girth. The conjecture predicts that even if that is not the case, we should be able to find spanning subgraphs of $G_p$ of large girth which are still good expanders. In this example we can do something slightly weaker. We can find a sequence of expanders $H_p$ on the same vertex set as $G_p$, of girth proportional to $\log p$, such that every pair of vertices in $H_p$ connected by an edge is of distance at most $O(1)$ in $G_p$. To find such graphs take a finite power of $S$ that contains two elements, say $a$ and $b$, that generate a free Zariski dense subgroup. This is always possible by the Tits alternative. Put $H_p:=\Cay(\SL(2,\mathbb Z/p\mathbb Z), \{a,b,a^{-1},b^{-1}\}).$ By \cite{BoGa} $H_p$ is an expander sequence. Since $\langle a,b\rangle$ is free, the girth must go to infinity. Indeed, any relation on $a,b$ that holds modulo infinitely many primes would have to hold in $\SL(2,\mathbb Z)$, which contradicts the freeness.

\item
Let's end by recalling an old conjecture: there is no sequence of finite bounded degrees graphs growing  in size  to infinity, so that all the induced balls in all the graphs in the sequence, are uniform expanders. 

For a related progress see \cite{FV}.

\end{enumerate}


\begin{thebibliography}{BKC}

 



\bibitem{abs}
N. Alon, I. Benjamini, and A. Stacey.
Percolation  on finite graph and isoperimetric ineqilities.
Ann. Prob. 32 (2004), 1727--1745.

\bibitem{pyond}
I. Benjamini, and O. Schramm,
Percolation beyond $\Z^d$ many questions and a few answers.
ECP 1 (1996), 71--82.

\bibitem{BS}
I. Benjamini, and O. Schramm,
Every graph with a positive Cheeger constant contains a tree with a positive Cheeger constant.
GAFA (1997), 403--419.

\bibitem{BLS}
I. Benjamini, R. Lyons, and O. Schramm,
Percolation perturbations in potential theory and random walks. Random walks and discrete potential theory (Cortona, 1997), 56–84, 
Sympos. Math., XXXIX, Cambridge Univ. Press, Cambridge, 1999. 

\bibitem{BoGa} 
J. Bourgain, and A. Gamburd, 
Uniform expansion bounds for Cayley graphs of $\SL_2(\mathbb F_p)$.
Ann of Math, 167, (2008), 625--642

\bibitem{BG}
E. Breuillard, and T. Gelander,
A topological Tits alternative.
Ann of Math, 166, (2007), 427--474


\bibitem{FV}
M. Fraczyk, and W. van Limbeek
Heat kernels are not uniform expanders.
arXiv:1905.13584


\bibitem{HLW}
S. Hoory, N. Linial, and A. Wigderson, Expander graphs and their applications. Bull. Amer. Math. Soc. (N.S.) 43 (2006), no. 4, 439–561. 

\bibitem{L}
A. Lubotzky, 
High dimensional expanders. Proceedings of the International Congress of Mathematicians—Rio de Janeiro 2018. Vol. I. Plenary lectures, 705–730, World Sci. Publ., Hackensack, NJ, 2018. 

\bibitem{PS}
I. Pak, and T. Smirnova-Nagnibeda, 
On non-uniqueness of percolation on nonamenable Cayley graphs. 
C. R. Acad. Sci. Paris Sér. I Math. 330 (2000), no. 6, 495


\end{thebibliography}
\end{document}